\magnification=\magstep1
\input amstex
\documentstyle{amsppt}

\define\defeq{\overset{\text{def}}\to=}
\define\ab{\operatorname{ab}}
\define\pr{\operatorname{pr}}
\define\tor{\operatorname{tor}}
\define\Gal{\operatorname{Gal}}
\define\Ker{\operatorname{Ker}}
\define\genus{\operatorname{genus}}

\define\Aut{\operatorname{Aut}}
\define\id{\operatorname{id}}
\define\Pic{\operatorname{Pic}}
\def \isom {\overset \sim \to \rightarrow}
\define\Sect{\operatorname{Sect}}
\define\Spec{\operatorname{Spec}}

\def\Card{\operatorname{Card}}
\define\res{\operatorname{res}}
\define\cor{\operatorname{cor}}

\define\Primes{\frak{Primes}}
\NoRunningHeads
\NoBlackBoxes
\topmatter

\title
Further examples of non-geometric sections of arithmetic fundamental groups
\endtitle

\author
Mohamed Sa\"\i di
\endauthor

\abstract
We show the existence of group-theoretic sections of certain {\it geometrically pro-nilpotent by abelian}
arithmetic fundamental groups of hyperbolic curves over $p$-adic local fields which are {\it non-geometric}, i.e.,
which do not arise from rational points. Among these quotients is the {\it geometrically metabelian} arithmetic fundamental group.
\endabstract

\endtopmatter

\document

\subhead
\S 0. Introduction/Statement of the Main Result
\endsubhead
Grothendieck's anabelian section conjecture predicts that sections of arithmetic fundamental groups of hyperbolic curves over finitely generated fields over $\Bbb Q$ arise from rational points 
(cf. [Sa\"\i di] for a more precise formulation of the conjecture). Accordingly, sections of arithmetic fundamental groups of hyperbolic curves over $p$-adic local fields; which are defined over number fields, and which arise from global sections, should arise from rational points. In this context it is tempting to predict 
a {\it $p$-adic analog of Grothendieck's anabelian section conjecture}. In [Sa\"\i di1] we investigated such analog, and exhibited two necessary and sufficient conditions for a section of the 
arithmetic fundamental group of a hyperbolic curve over a $p$-adic local field to be geometric, i.e., to arise from a rational point (cf. loc. cit. Theorem 4.5).

For the time being there are no examples of sections of (the full) arithmetic fundamental groups of hyperbolic curves over $p$-adic local fields which are non-geometric, and one can still hope the validity of a $p$-adic analog of the section conjecture.
On the other hand, recent examples were found of group-theoretic sections of certain {\it (geometrically characteristic) quotients} of 
arithmetic fundamental groups of curves over $p$-adic local fields which are {\it non-geometric}. Hoshi constructed examples of sections of the {\it geometrically pro-$p$} quotient
of arithmetic fundamental groups of curves over $p$-adic local fields which are non-geometric (cf. [Hoshi]).
(Actually, Hoshi's example arises from group-theoretic sections of geometrically pro-$p$ fundamental groups
of hyperbolic curves over number fields (cf. loc. cit.).) In [Sa\"\i di]  we constructed examples of group-theoretic sections of 
{\it geometrically prime-to-$p$} fundamental groups
of hyperbolic curves over $p$-adic local fields which are non-geometric (cf. loc. cit. $\S3$). Further, in
[Sa\"\i di2] we provided examples of group-theoretic sections of 
the "{\it \'etale by geometrically abelian}" fundamental group
of hyperbolic curves over $p$-adic local fields which are non-geometric.
 {\it The existence of these examples is crucial for our understanding 
of the $p$-adic section conjecture}. Indeed, if the $p$-adic version of the section conjecture
holds true then it may possibly hold true even for smaller quotients of the arithmetic fundamental group,
and one would like to know these quotients in this case. On the other hand, more elaborate examples
of non-geometric sections as above may lead to a counterexample for the $p$-adic version of the section conjecture.

In this note we provide further examples of sections of certain quotients of arithmetic fundamental groups of 
curves over $p$-adic local fields which are non-geometric. 
These quotients include the {\it geometrically metabelian} and certain {\it geometrically pro-nilpotent by abelian} quotients.

Next, we fix notations and state our main results.

\smallskip
$\bullet$\ Let 
$$1\to H'\to H @>{\pr}>> G \to 1$$ 
be an exact sequence of profinite groups. We will refer to
a continuous homomorphism $s:G\to H$ satisfying $\pr \circ s=\id_{G}$ as a (group-theoretic)
{\bf section}, or {\bf splitting}, of the above sequence, or simply a section of the projection $\pr : H\twoheadrightarrow G$.
We denote by $\Sect\left(H\twoheadrightarrow G\right)$ the set of sections of the projection $H\twoheadrightarrow G$.

\smallskip
$\bullet$ Given a profinite group $H$, and a prime integer $\ell$, we will denote by $H^{\ell}$ the maximal {\it pro-$\ell$} quotient of $H$,
$H^{\ab}$ the maximal {\it abelian} quotient of $H$, and
$H^{\ab,\ell}$ its maximal {\it abelian pro-$\ell$} quotient. Thus $H^{\ab,\ell}=(H^{\ell})^{\ab}$.

\smallskip
Let $p\ge 2$ be a {\it prime integer}, 
and $k$ a $p$-{\bf adic local field}; meaning $k/\Bbb Q_p$ is a finite extension,
with ring of integers $\Cal O_k$, and residue field $F$. 
Thus $F$ is a {\it finite field} of characteristic $p$. 
Let $X\to \Spec k$ 
be a proper, smooth, and geometrically connected {\bf hyperbolic} (i.e., $\genus (X)\ge 2$) {\bf curve} over $k$.
Let $\eta$ be a geometric point of $X$ above its generic point, which
determines an algebraic closure $\overline k$
of $k$, and a geometric point $\bar {\eta}$ of $\overline {X} \defeq X\times _k \overline k$.
There exists a canonical exact sequence of profinite groups (cf. [Grothendieck], Expos\'e IX, Th\'eor\`eme 6.1)
$$1\to \pi_1(\overline {X},\bar \eta)\to \pi_1(X, \eta) \to G_k\to 1.$$
Here, $\pi_1(X, \eta)$ denotes the {\it arithmetic \'etale fundamental group} of $X$ with base
point $\eta$, $\pi_1(\overline {X},\bar \eta)$ the \'etale fundamental group of $\overline {X}\defeq 
X\times _k \overline k$ with base
point $\bar \eta$, and $G_k\defeq \Gal (\overline k/k)$ the absolute Galois group of $k$.

\smallskip
$\bullet$ Let $\Pi$ be a quotient of $\pi_1(X,\eta)$ such that the projection $\pi_1(X, \eta) \twoheadrightarrow G_k$ factors as $\pi_1(X,\eta)\twoheadrightarrow \Pi \twoheadrightarrow G_k$, and which is geometrically non-trivial; meaning $\Ker (\Pi\twoheadrightarrow G_k)$ is non-trivial.
Given a section $s:G_k\to \Pi$ of the projection $\Pi\twoheadrightarrow G_k$, we say that $s$ is {\bf geometric} if $s(G_k)$ is contained in (hence equal to) the decomposition 
group $D_x\subset \Pi$ associated to a rational point $x\in X(k)$. In this case we say $s$ arises from the rational point $x$.
We say that the section $s$ is {\bf non-geometric} if $s$ is not geometric in the above sense, i.e., $s(G_k)$ is {\it not}
contained in the decomposition group associated to a rational point $x\in X(k)$. (Note that in the above discussion the decomposition group $D_x$ is only defined up to conjugation by elements
of $\Ker (\Pi\twoheadrightarrow G_k)$.)

\smallskip
$\bullet$ 
We {\bf assume} $X(k)\neq \emptyset$. We {\bf fix} a $k$-rational point $x\in X(k)$, and $s\defeq s_x:G_k\to \pi_1(X,\eta)$ a section of the projection $\pi_1(X,\eta)\twoheadrightarrow G_k$ associated to $x$. Thus $s$ is defined only up to conjugation by $\pi_1(\overline {X},\bar \eta)$. 
Note that the section $s$ induces a structure of $G_k$-group on 
any characteristic quotient of $\pi_1(\overline {X},\bar \eta)$.

\smallskip 
$\bullet$ Let $\Delta$ be a quotient  of $\pi_1(\overline {X},\bar \eta)$ which fits in an exact sequence
$$1\to \widetilde H\to \pi_1(\overline {X},\bar \eta)\to \Delta\to 1,$$ 
where $\widetilde H\defeq \Ker (\pi_1(\overline {X},\bar \eta)\twoheadrightarrow \Delta)$. We consider the following condition $(\star)$ on $\Delta$.

\proclaim {Condition ($\star$)} 

\smallskip
{\bf (i)} $\Delta$ is {\bf pro-nilpotent}, and is a {\bf characteristic} quotient of $\pi_1(\overline {X},\bar \eta)$.

\smallskip
{\bf (ii)} $H^0(U,\Delta)=0$ for every open subgroup $U$ of $G_k$.

\smallskip
{\bf (iii)} The quotient  $\pi_1(\overline {X},\bar \eta)\twoheadrightarrow \pi_1(\overline {X},\bar \eta)^{\ab}$ {\bf factors} as 
$$\pi_1(\overline {X},\bar \eta)\twoheadrightarrow \Delta \twoheadrightarrow \pi_1(\overline {X},\bar \eta)^{\ab}.$$

\smallskip
{\bf (iv)}\ Let $H\defeq \widetilde H^{\ab}$, and $\Gamma\defeq \pi_1(\overline X,\bar \eta)/\Ker (\widetilde H\twoheadrightarrow H)$. We have a push-out diagram
$$
\CD
1 @>>> \widetilde H @>>> \pi_1(\overline X,\bar \eta) @>>> \Delta @>>> 1\\
@. @VVV @VVV @|\\
1 @>>> H @>>> \Gamma @>>> \Delta @>>> 1\\
\endCD
$$ 
where the middle and left vertical maps are surjective. There {\bf exists a prime integer $\ell\neq p$}, such that the natural surjective map
$$\Gamma^{\ell}\twoheadrightarrow \Delta^{\ell}$$ 
is {\bf not an isomorphism}.
 \endproclaim

\smallskip
Let $\Pi\defeq \pi_1(X,\eta)/\Ker (\pi_1(\overline X,\bar \eta)\twoheadrightarrow \Gamma)$, and  
$\widetilde \Pi\defeq \pi_1(X,\eta)/\Ker (\pi_1(\overline X,\bar \eta)\twoheadrightarrow \Delta)$.
We have the following push-out diagrams
  $$
 \CD
 1@>>>   \pi_1(\overline {X},\bar \eta)  @>>> \pi_1(X,\eta) @>>> G_k @>>> 1 \\
 @. @VVV @VVV   @|   \\
 1@>>> \Gamma @>>> \Pi @>>> G_k @>>> 1 \\
 @. @VVV @VVV  @|   \\
 1@>>> \Delta @>>> \widetilde \Pi @>>> G_k @>>> 1 \\
 \endCD
 $$
where the vertical maps are surjective. Thus $\Ker (\Pi\twoheadrightarrow \widetilde \Pi)=\Ker (\Gamma \twoheadrightarrow \Delta)=H$.

 \smallskip
$\bullet$ Examples of quotients $\Delta$ satisfying the condition $(\star)$ are: $\Delta=\pi_1(\overline {X},\bar \eta)^{\ab}$ 
is the maximal abelian quotient of $\pi_1(\overline {X},\bar \eta)$ (cf. [Sa\"\i di3], Lemma 1.3, for the condition $(\star)$(ii)); in this case 
$\widetilde \Pi\defeq \pi_1(X,\eta)^{(\ab)}$ is the geometrically abelian quotient of $\pi_1(X,\eta)$, $\Gamma$ is the maximal metabelian quotient of $\pi_1(\overline {X},\bar \eta)$,
and $\Pi$ is the {\bf geometrically metabelian} quotient of $\pi_1(X,\eta)$. More generally, any pro-nilpotent characteristic quotient $\Delta$ of $\pi_1(\overline {X},\bar \eta)$ 
which satisfies conditions $(\star)$(ii), $(\star)$(iii),
and for which there exists a prime integer $\ell\neq p$ such that the natural projection 
$\pi_1(\overline {X},\bar \eta)^{\ell}\twoheadrightarrow \Delta ^{\ell}$ is not an isomorphism,
satisfies the condition $(\star)$.

Given a finite extension $k'/k$ (all finite extensions of $k$ we consider are contained in $\overline k$), and the corresponding open subgroup $G_{k'}\subseteq G_k$, 
we will denote by $\Pi_{k'}$ the pullback of the group extension $\Pi$ by $G_{k'}\hookrightarrow G_k$.
Thus we have a commutative diagram of exact sequences
$$
\CD
1@>>> \Gamma @>>> \Pi_{k'} @>>> G_{k'} @>>> 1 \\
@. @| @VVV @VVV\\
1@>>> \Gamma @>>> \Pi @>>> G_k @>>> 1 \\
\endCD
$$
where the right square is cartesian. Likewise we write $\widetilde \Pi_{k'}\defeq \widetilde \Pi\times _{G_k}G_{k'}$.
Note that $\Pi_{k'}$ is a quotient of $\pi_1(X_{k'},\eta)$, where $X_{k'}\defeq X\times _{\Spec k}\Spec k'$ and $\eta$ is naturally induced by the above geometric point $\eta$.
Our first main result in this paper is the following.

\proclaim {Theorem A} We use notations as above. Let $X$ be a proper, smooth, and geometrically connected hyperbolic curve over the $p$-adic local field $k$.
Assume that $X(k)\neq \emptyset$, and $X$ has {\bf potentially good reduction}.
Let $\Delta$ be a quotient of $\pi_1(\overline {X},\bar \eta)$ {\bf satisfying the condition} $(\star)$, and $\Pi$ the corresponding quotient of $\pi_1(X,\eta)$ as in the above
discussion which fits in the exact sequence $1\to \Gamma \to \Pi \to G_k\to 1$.
Then there exists a finite extension $\tilde k/k$ such that the following holds. For every finite extension $k'/\tilde k$, there exists  
a section $s:G_{k'}\to \Pi_{k'}$ of the projection $\Pi_{k'}\twoheadrightarrow G_{k'}$ which is {\bf non-geometric}.
\endproclaim

As a corollary of Theorem A we obtain the following (cf. examples discussed after introducing the condition $(\star)$).

\proclaim {Corollary B} There exist {\bf non-geometric} sections of geometrically metabelian arithmetic fundamental groups of hyperbolic curves over $p$-adic local fields. 
\endproclaim

Let $m\ge 1$ be an integer. With the notations above, let $\Delta_m\defeq \Delta_{m,X}$ be the maximal {\bf $m$-step solvable pro-$p$} quotient of $\pi_1(\overline {X},\bar \eta)$, and 
$\Pi_m\defeq \Pi_{m,X}$ the {\bf geometrically $m$-step solvable pro-$p$} quotient of $\pi_1(X,\eta)$ which sits in the exact sequence 
$$1\to \Delta_m\to \Pi_m\to G_k\to 1$$ 
(cf. [Sa\"\i di3], $\S1$). Note that $\Delta _m$ doesn't satisfy condition $(\star)$(iii). 
It is plausible, in light of Hoshi's example in [Hoshi] (cf. above discussion), that there exist non-geometric sections of $\Pi_m$ for a suitable $X/k$ as above  (this is easily seen if $m=1$, using the Kummer exact sequence associated to $\Pic^0_X$).
In this context we prove the following.

\proclaim {Theorem C} We use notations as above. There  exists an integer $N\ge 2$, such that the following holds. For every prime integer $p\ge N$ there exists 
 a proper, smooth, and geometrically connected hyperbolic curve $X$ over a $p$-adic local field $k$, an integer $m\ge 2$, and a section $s:G_k\to \Pi_{m,X}$ of the projection $\Pi_{m,X}\twoheadrightarrow G_k$ which is {\bf non-geometric}.
\endproclaim

\subhead
Acknowledgment
\endsubhead
I would like to thank Akio Tamagawa for very fruitful discussions during which the proofs of Theorem A and Theorem C were established, while the author was visiting the Research Institute for Mathematical Sciences at Kyoto university (RIMS), the idea of proof of Theorem C is due to him. I would also like to thank the referee for his/her very helpful comments.

\subhead
\S 1. Proof of Theorem A
\endsubhead
We use the notations introduced in $\S0$, as well as the notations and assumptions in Theorem A.
Thus $X$ is a proper, smooth, and geometrically connected hyperbolic curve over the $p$-adic local field $k$, $X(k)\neq \emptyset$, and we assume (without loss of generality) that 
$X$ has good reduction over $\Cal O_k$. Further,
$\Delta$ is a characteristic quotient of $\pi_1(\overline {X},\bar \eta)$ satisfying the condition $(\star)$, and $\Pi$ is the corresponding quotient of $\pi_1(X,\eta)$ as above
which fits in the exact sequence $1\to \Gamma \to \Pi \to G_k\to 1$ (cf. $\S0$). We have the following commutative diagram of exact sequences.

$$
 \CD
 @. 1 @. 1\\
 @. @VVV @VVV\\
 @.  H  @= H @. \\
 @. @VVV @VVV \\
 1@>>> \Gamma @>>> \Pi @>>> G_k @>>> 1 \\
 @. @VVV @VVV @||\\
 1@>>> \Delta @>>> \widetilde \Pi @>>> G_k @>>> 1 \\
 @. @VVV @VVV \\
@. 1 @. 1 \\
 \endCD
 \tag 1.1
 $$

Recall $x\in X(k)$ is a $k$-rational point, and $s\defeq s_x:G_k\to \pi_1(X,\eta)$ a section of the projection $\pi_1(X,\eta)\twoheadrightarrow G_k$ associated to $x$.
Further, $s$ induces sections: $s_1\defeq s_{1,x}:G_k\to \widetilde \Pi$ of the projection
$\widetilde \Pi\twoheadrightarrow G_k$, and $s_2\defeq s_{2,x}:G_k\to \Pi$ of the projection
$\Pi\twoheadrightarrow G_k$, which fit in a commutative diagram
$$
\CD
G_k @>s_2>> \Pi \\
@|  @VVV \\
G_k @>s_1>> \widetilde \Pi \\
\endCD
$$ 
where the right vertical map is the one in diagram (1.1).

The profinite group $\Delta $ is {\bf finitely generated}, as follows from the well-known finite generation of the profinite group 
$\pi_1(\overline X,\overline \eta)$ which projects onto $\Delta$.
Let $\{\Delta ^{i}\}_{i\ge 1}$ be a countable system of {\bf characteristic open} subgroups of $\Delta$ such that 
$$\Delta ^{i+1}\subseteq \Delta ^{i},\ \ \ \ \Delta ^{1} = \Delta, \ \ \ \ \text {and} \ \ \bigcap _{i\ge 1}\Delta ^{i} = \{1\}.$$
Write $\Delta _i\defeq \Delta /\Delta ^{i}$. Thus $\Delta _i$ is a {\it finite characteristic} quotient of $\Delta $, and
we have a push-out diagram of exact sequences

$$
\CD
1@>>> \Delta @>>>   \widetilde \Pi  @>>>  G_k @>>> 1\\
@.  @VVV    @VVV    @| \\
1@>>>  \Delta _i  @>>>  \Pi _i @>>>  G_k @>>> 1\\
\endCD
\tag 1.2
$$
which defines a ({\it geometrically finite}) quotient $\Pi_i$ of $\widetilde \Pi$.
The section $s_1$  induces a section
$$\rho_i:G_k\to \Pi_i$$
of the projection $\Pi_i\twoheadrightarrow G_k$, $\forall i\ge 1$. 
Write 
$$\widetilde \Pi^i\defeq \widetilde \Pi^i[s_1] \defeq \Delta ^{i}.s_1 (G_k).$$
Note that $\widetilde \Pi^i\subseteq \widetilde \Pi$ is an open subgroup which contains the image 
$s_1 (G_k)$ of $s_1$. Write $\Pi^i$ for the inverse image of $\widetilde \Pi ^i$
in $\pi_1(X,\eta)$. Thus $\Pi^i\subseteq \pi_1(X,\eta)$ is an open subgroup corresponding
to an \'etale cover 
$$X_i\to X_1\defeq X$$
defined over $k$ (since $\Pi^i$ maps onto $G_k$ via the natural projection $\pi_1(X,\eta)\twoheadrightarrow G_k$, by the very 
definition of $\Pi^i$). 

Note that the \'etale cover $\overline X_i\defeq X_i\times _{\Spec k}\Spec {\overline k}\to \overline X$ is Galois with
Galois group $\Delta _i$, and we have a commutative diagram of \'etale covers

$$
\CD
\overline X_i  @>>>   \overline X \\
@VVV    @VVV \\
X_i  @>>>  X\\
\endCD
$$
where $\overline X_i\to X$ is Galois with Galois group $\Pi_i$, and $\overline X_i\to X_i$ is Galois with Galois group $\rho_i(G_k)$. 
We have a commutative diagram of exact sequences
$$
\CD
@.   1@.    1@.   \\
@. @VVV  @VVV\\
1 @>>>  \widetilde \Delta ^{i}= \pi_1(\overline {X_i},\bar \eta)@>>> \Pi^i=\pi_1(X_i,\eta) @>>>  G_k @>>> 1\\
@.   @VVV    @VVV   @| \\
1 @>>> \pi_1(\overline {X},\bar \eta) @>>>   \pi_1(X,\eta) @>>>  G_k  @>>> 1 \\
\endCD
$$
where $\widetilde \Delta^{i}$ is the inverse image of $\Delta ^{i}$
in $\pi_1(\overline {X},\bar \eta)$, and the equalities $\widetilde \Delta ^{i}= \pi_1(\overline {X_i},\bar \eta)$, $\Pi^i=\pi_1(X_i,\eta)$,
are natural identifications; the base points $\eta$ (resp. $\overline \eta$) of $X_i$ (resp. $\overline X_i$) are those induced by 
the base points $\eta$ (resp. $\overline \eta$) of $X$ (resp. $\overline X$). Note that $\Pi^{i+1}\subseteq \Pi^i$, and $\widetilde \Delta ^{i+1}\subseteq \widetilde \Delta ^i$, 
as follows from the various definitions.

\proclaim {Lemma 1.1} With the above notations and those in $\S0$, the following holds:
$$\widetilde H=\bigcap _{i\ge 1} \widetilde \Delta ^i.$$
\endproclaim

\demo{Proof}
Follows from the various definitions.
\qed
\enddemo

We take this opportunity to correct a mistake that occurred in [Sa\"\i di2], Lemma 1.1. The claim there that $\Cal I_X=\bigcap_{i\ge 1} \Pi_i$ is false, however this doesn't affect the validity of the results or other assertions made in loc. cit..

For each integer $i\ge 1$, consider the push-out diagram
$$
\CD
1 @>>>  \widetilde \Delta ^i= \pi_1(\overline {X_i},\bar \eta)@>>> \Pi^i=\pi_1(X_i,\eta) @>>>  G_k @>>> 1\\
@.   @VVV    @VVV   @|\\
1 @>>>  \widetilde \Delta ^{i,\ab}@>>> \Pi^{(i,\ab)}@>>>  G_k @>>> 1\\
\endCD
$$
where $\widetilde \Delta ^{i,\ab}$ is the maximal abelian quotient of $\widetilde \Delta ^{i}$, and
$\Pi^{(i,\ab)}$ is the {\bf geometrically abelian} fundamental group of $X_i$.
Consider the commutative diagram 

$$
\CD
1 @>>> H  @>>> \Cal H \defeq   \Cal H [s_1] @>>>  G_k @>>> 1\\
@.   @VVV    @VVV   @V{s_1}VV\\
1 @>>>   H  @>>> \Pi @>>>  \widetilde \Pi  @>>> 1 \\ 
\endCD
\tag 1.3$$ 
 where the right square is cartesian. Thus (the group extension) $\Cal H$ is the pull-back of (the group extension) 
 $\Pi$ via the section  $s_1:G_k\to \widetilde \Pi$.

 \proclaim {Lemma 1.2} We have natural identifications $H  \isom \underset{i\ge 1} \to{\varprojlim}\ \widetilde \Delta ^{i,\ab}$, and
 $\Cal H \isom \underset{i\ge 1} \to{\varprojlim}\  \Pi^{(i,\ab)}$.
 \endproclaim
 
 \demo {Proof} Similar to the proof of Lemma 1.3 in [Sa\"\i di2].
 \qed
 \enddemo

The section $s_2:G_k\to \Pi$, which lifts the section $s_1$, induces a section $s_2:G_k\to \Cal H$ of the projection $\Cal H\twoheadrightarrow G_k$ (since $s_2(G_k)\subset \Cal H$).
We fix the section $s_2:G_k\to \Cal H$ as a  base point of the torsor of splittings of the upper sequence in diagram (1.3). Thus the set of splittings of the group extension $\Cal H$, modulo conjugation by elements of $H$, is a torsor under $H^1(G_k,H)$; the $G_k$-module structure of $H$ being deduced from diagram (1.3). 
The splitting $s_2:G_k\to \Cal H$ thus corresponds to $0\in H^1(G_k,H)$. Note that the set of splittings $\tilde s:G_k\to \Cal H$ of the group extension $\Cal H$ is in one-to-one correspondence with the set of sections $\tilde s:G_k\to \Pi$ 
of the projection $\Pi\twoheadrightarrow G_k$ which lift the section $s_1$.

Let $\tilde s:G_k\to \Cal H$ be a section of the group extension $\Cal H$, which induces a  section $\tilde s:G_k\to \Pi$ 
of the projection $\Pi\twoheadrightarrow G_k$ which lift the section $s_1$. Let $[\tilde s]$ be the class of $\tilde s$ in $H^1(G_k,H)$ (cf. above discussion).

\proclaim {Fact 1.3} The section $\tilde s:G_k\to \Pi$ is geometric if and only if $[\tilde s]=0$.  In this case the section $\tilde s$ is associated to the rational point $x\in X(k)$.
 \endproclaim
 
 \demo {Proof} First, assume that the section $\tilde s:G_k\to \Pi$ is geometric and arises from a rational point $\tilde x\in X(k)$. Both sections $\tilde s:G_k\to \Pi$, and 
 $s_2:G_k\to \Pi$, induce splittings $\tilde s^{\ab}:G_k\to \pi_1(X,\eta)^{(\ab)}$, and  $s_2^{\ab}:G_k\to \pi_1(X,\eta)^{(\ab)}$ of the group extension $1\to \pi_1(\overline {X},\bar \eta)^{\ab}  \to \pi_1(X,\eta)^{(\ab)} \to G_k \to 1$,
 where  $\pi_1(X,\eta)^{(\ab)}$ is the geometrically abelian quotient of  $\pi_1(X,\eta)$. Further one has   $\tilde s^{\ab}=s_2^{\ab}$ (cf. condition $(\star)$(iii), and the fact that both $\tilde s$ and $s_2$ lift the section $s_1$). A standard argument, resorting to the Kummer exact sequence associated to the jacobian $\Pic^0_X$ of $X$, shows that $\tilde x=x$ (cf. [Tamagawa], Proposition 2.8). 

Next, we claim $[\tilde s]=0$. Indeed, the classes of $\tilde s$ and $s_2$ in $H^1(G_k,\Gamma)$ coincide as both sections are geometric and associated to the same rational point $x$, hence 
 $\tilde s$ and $s_2$  are conjugate by an element of $\Gamma$. Here we view the set of splittings of the group extension $\Pi$ (of $\Gamma$ by $G_k$) as a torsor under $H^1(G_k,\Gamma)$, with base point the class
 of the section $s_2$. Further the natural map $H^1(G_k,H)\to H^1(G_k,\Gamma)$ of pointed cohomology sets is injective as follows from the 
 condition $(\star)$(ii) (cf. [Serre], I.$\S5$, Proposition 38, and diagram (1.1)). (Here the $G_k$-module structure on $H$ (resp. $G_k$-group structure on $\Gamma$) is induced by the section $s_1$ (resp. $s_2$) (cf. diagram (1.1)).) Thus $[\tilde s]=0$.

 Conversely, if $[\tilde s]=0$, then $\tilde s$ is conjugate to $s_2$ by an element of $H$ hence is geometric and associated to the rational point $x$.
 \qed
 \enddemo

As a consequence we obtain the following.

\proclaim {Lemma 1.4} Let $k'/k$ be a finite extension. There exists a  section $\tilde s:G_{k'}\to \Pi_{k'}$ of the projection $\Pi_{k'}\twoheadrightarrow G_{k'}$, which lifts the section
$s_{1,k'}:G_{k'}\to \widetilde \Pi_{k'}$ of  the projection $\widetilde \Pi_{k'}\twoheadrightarrow G_k$ induced by $s_1$, and
 which is {\bf non-geometric}, if and only if $H^1(G_{k'},H)\neq 0$. 
 \endproclaim
 
Thus proving Theorem A reduces to proving the following.

\proclaim {Proposition 1.5} There exists a finite extension $\tilde k/k$ such that $H^1(G_{k'},H)\neq 0$ for every finite extension $k'/\tilde k$. 
 \endproclaim

The rest of this section is devoted to proving Proposition 1.5. Let $\ell\neq p$ be a prime integer such that 
the map $\Gamma^{\ell}\twoheadrightarrow \Delta^{\ell}$ is not an isomorphism (cf. condition $(\star)$(iv)).
Write $\pi_1(\overline {X},\bar \eta)^{\ell}$ for the maximal {\it pro-$\ell$} quotient of $\pi_1(\overline {X},\bar \eta)$, and 
 $$\pi_1(X, \eta)^{(\ell)}\defeq \pi_1(X,\eta)/\Ker  (\pi_1(\overline {X},\bar \eta) \twoheadrightarrow \pi_1(\overline {X},\bar \eta)^{\ell})$$ 
 for the {\it geometrically pro-$\ell$} quotient of 
 $\pi_1(X,\eta)$, which fits in the exact sequence
 $$1\to \pi_1(\overline {X},\bar \eta)^{\ell}\to \pi_1(X, \eta)^{(\ell)} \to G_k\to 1.$$

Let $s^{\ell}=s^{\ell}_x:G_k\to \pi_1(X,\eta)^{(\ell)}$ be the section of the projection $\pi_1(X, \eta)^{(\ell)} \twoheadrightarrow G_k$ induced by the section $s=s_x$.
This section induces a representation 
$$\rho^{\ell} : G_k\to \Aut (\pi_1(\overline {X},\bar \eta)^{\ell})$$ 
which factors as $G_k\twoheadrightarrow G_F \to \Aut (\pi_1(\overline {X},\bar \eta)^{\ell})$,
 where $G_F$ is the quotient of $G_k$ by its inertia subgroup, since $X$ has good reduction over $\Cal O_k$.
 Further the image of the representation $\rho^{\ell}$ is almost pro-$\ell$, i.e., $\rho^{\ell}(G_k)$ possesses an open subgroup which is pro-$\ell$. In particular, there exists a finite extension
 $\tilde k/k$ such that the restriction $\rho^{\ell}_{\tilde k}:G_{\tilde k}\to \Aut (\pi_1(\overline {X},\bar \eta)^{\ell})$ of $\rho^{\ell}$ to $G_{\tilde k}$ has a pro-$\ell$ image. In order to prove 
 Proposition 1.5 we will, without loss of generality, {\bf assume that the image of $\rho^{\ell}$ is pro-$\ell$}, and will show  $H^1(G_{k},H)\neq 0$.

 Let $\Delta ^{\ell}$ be the maximal pro-$\ell$ quotient of $\Delta$, and $\widetilde \Pi^{(\ell)}\defeq \widetilde \Pi /\Ker (\Delta \twoheadrightarrow \Delta ^{\ell} )$
 the geometrically pro-$\ell$ quotient of $\widetilde \Pi$. We have the following commutative diagram of exact sequences
 $$
 \CD
 1@>>> \Delta @>>>   \widetilde \Pi  @>>>  G_k @>>> 1\\
 @. @VVV  @VVV @|\\
 1@>>> \Delta^{\ell} @>>>   \widetilde \Pi ^{(\ell)}  @>>>  G_k @>>> 1\\
 \endCD
 $$
 where the left and middle vertical maps are surjective.
 For $i\ge 1$, let $N ^i$ be the image of $\Delta ^i$ in $\Delta ^{\ell}$, and $\widetilde N^i$, $\widehat N^i$,
 the pre-images of $N^i$ in $\pi_1(\overline {X},\bar \eta)$, and $\pi_1(\overline {X},\bar \eta)^{\ell}$; respectively. Note that  
$\widetilde N^i$ is a characteristic subgroup of $\pi_1(\overline {X},\bar \eta)$, and $\widehat N^i$ is stable by the action of $s^{\ell}(G_k)$.

Let $\widetilde U^i\defeq \widetilde N^i.s(G_k)$, and $\widehat U^i\defeq \widehat N^i.s^{\ell}(G_k)$, for $i\ge 1$.
Thus $\widetilde U^i$ is an open subgroup of $\pi_1({X},\eta)$ 
corresponding to an \'etale cover $Y_i\to X$, and the \'etale cover $X_i\to X$ factorises as 
 $$X_i\to Y_i\to X$$
 [$\pi_1(X_i,\eta)\subset \pi_1(Y_i,\eta)$ as follows from the various definitions], where $X_i\to Y_i$ is an \'etale cover of degree prime-to-$\ell$, since $\Delta$ 
 (hence also $\Delta_i$, for $i\ge 1$)
 is pro-nilpotent (see condition $(\star)$(i)). 
 Further $\widetilde U^i$, and $\widehat U^i$, are naturally identified with $\pi_1(Y_i,\eta)$, and $\pi_1(Y_i,\eta)^{(\ell)}$; respectively, $\forall i\ge 1$. Here $\pi_1(Y_i,\eta)^{(\ell)}$ is the geometrically pro-$\ell$ quotient of 
 $\pi_1(Y_i,\eta)$, and sits in an exact sequence 
 $$1\to \pi_1(\overline {Y}_i,\bar \eta)^{\ell}\to \pi_1(Y_i,\eta)^{(\ell)} \to G_k\to 1,$$ 
 where $\overline Y_i\defeq Y\times _k\overline k$.

 The natural action of $G_k$ on $\pi_1(\overline X,\bar \eta)^{\ell}$, and which factorises through $G_F^{\ell}$ by our assumption on the representation $\rho^{\ell}$,
 is compatible with its action on the open subgroup $\pi_1(\overline Y_i,\bar \eta)^{\ell}$, hence this latter action also factorises through $G_F^{\ell}$.

 There is a surjective homomorphism (recall Lemma 1.2)
$$H^1(G_k,H)=\varprojlim _{i\ge 1} H^1(G_k, \widetilde \Delta ^{i,\ab}) \twoheadrightarrow \varprojlim _{i\ge 1} H^1(G_k, \widetilde \Delta ^{i,\ab,\ell}).$$
(Indeed, $\varprojlim _{i\ge 1} H^1(G_k, \widetilde \Delta ^{i,\ab}) \isom \prod_{l\in \Primes} \varprojlim _{i\ge 1} H^1(G_k, \widetilde \Delta ^{i,\ab,l})$, where the product is over all prime integers $l$ and the above homomorphism is the projection onto the $\ell$-th factor.)
Further the \'etale covers $\{X_i\to Y_i\}_{i\ge 1}$ induce a homomorphism (recall $\widetilde \Delta ^{i}=\pi_1(\overline X_i,\bar \eta)$)
$$\varprojlim _{i\ge 1} H^1(G_k, \widetilde \Delta ^{i,\ab,\ell})\to \varprojlim _{i\ge 1} H^1(G_k, \pi_1(\overline Y_i,\bar \eta)^{\ab,\ell}),$$
which is surjective. More precisely, the map  $H^1(G_k, \widetilde \Delta ^{i,\ab,\ell})\to H^1(G_k, \pi_1(\overline Y_i,\bar \eta)^{\ab,\ell})$ is surjective, $\forall i\ge 1$,
as follows easily from a restriction-corestriction argument using the fact that the degree of the cover $X_i\to Y_i$ is prime-to-$\ell$ (observe the maps on cohomology induced by 
the natural maps $\pi_1(\overline Y_i,\bar \eta)^{\ab,\ell} @>\res>> \widetilde \Delta ^{i,\ab,\ell}@>\cor>> \pi_1(\overline Y_i,\bar \eta)^{\ab,\ell}$ arising from the morphisms
$\Pic^0(Y_i)\to \Pic^0(X_i)\to \Pic^0(Y_i)$, where the first one is the pull-back map of line bundles and the second is the norm map).
Thus in order to prove Proposition 1.5 it suffices to show the following. 

\proclaim {Proposition 1.6} With the above notations, it holds that 
 $$\varprojlim _{i\ge 1} H^1(G_k, \pi_1(\overline Y_i,\bar \eta)^{\ab,\ell})\neq 0.$$
\endproclaim

\demo{Proof} As discussed above the natural action of $G_k$ on $\pi_1(\overline Y_i,\bar \eta)^{\ab,\ell}$ factors through $G_F^{\ell}$ (which is isomorphic to $\Bbb Z_{\ell}$).
There is an injective inflation map 
$$\inf:\varprojlim _{i\ge 1} H^1(G_F^{\ell}, \pi_1(\overline Y_i,\bar \eta)^{\ab,\ell})\hookrightarrow \varprojlim _{i\ge 1} H^1(G_k, \pi_1(\overline Y_i,\bar \eta)^{\ab,\ell}).$$
Further  
$$(\varprojlim _{i\ge 1} \pi_1(\overline Y_i,\bar \eta)^{\ab,\ell})_{G_F^{\ell}}=H^1(G_F^{\ell},\varprojlim _{i\ge 1} \pi_1(\overline Y_i,\bar \eta)^{\ab,\ell}) =
 \varprojlim _{i\ge 1} H^1(G_F^{\ell}, \pi_1(\overline Y_i,\bar \eta)^{\ab,\ell})$$
where the notation $(\ \ )_{G_F^{\ell}}$ stands for the co-invariant module, the first equality follows from the fact that $G_F^{\ell}$ is procyclic, and the second
follows from [Neukirch-Schmidt-Winberg] (2.3.5) Corollary.

There is a natural isomorphism 
$$\left(\Ker \left(\Gamma ^{\ell}\to \Delta ^{\ell} \right) \right) \isom \varprojlim _{i\ge 1} \pi_1(\overline Y_i,\bar \eta)^{\ab,\ell}.$$ 
(Proof similar to the proof of Lemma 1.2.) Thus $\varprojlim _{i\ge 1} \pi_1(\overline Y_i,\eta)^{\ab,\ell}\neq 0$ (cf. condition $(\star)$(iii), and our choice of $\ell$). 
The proof of Proposition 1.6 follows from the following.

\proclaim {Lemma 1.7} Let $T$ be an abelian pro-$\ell$ group, and $P$ an infinite pro-$\ell$ cyclic group. Assume $T$ is a continuous $P$-module. 
Then the co-invariant module $(T)_P=\{0\}$ is trivial if and only if $T=\{0\}$ itself is trivial. 
\endproclaim

\demo{Proof} Let $T^{\wedge}$ be the Pontryagin dual of $T$ which is an $\ell$-primary torsion group. The dual of $(T)_P$ is the invariant group $(T^{\wedge})^P$.
It suffices to show that $(T^{\wedge})^P$ is trivial if and only if $(T^{\wedge})$ is trivial.
The action of $P$ on $T^{\wedge}$ is discrete, in particular $T^{\wedge}$ is the union of finite $\ell$-groups which are stable $P$-submodules. We can thus reduce to the case where 
$T^{\wedge}$
and $P$ are finite. Suppose $(T^{\wedge})$ is finite, and non-trivial, then $(T^{\wedge})^P$ is non-trivial since its order is divisible by $\ell$, and $(T^{\wedge})^P$ contains $0$.
\qed
\enddemo

This finishes the proof of Proposition 1.6, hence the proof of Proposition 1.5, and
the proof of Theorem A.
\qed
\enddemo

\subhead
\S 2. Proof of Theorem C
\endsubhead
The rest of this paper is devoted to proving Theorem C. We use the notations as introduced in $\S0$, and the statement of Theorem C.

Let $K$ be a number field (finite extension of $\Bbb Q$), and $\overline K$ an algebraic closure of $K$. Let
$X$ be a proper, smooth, and geometrically connected hyperbolic curve over $K$. Write
$J\defeq \Pic^0_X$ for the jacobian of $X$. Assume $X(K)\neq \emptyset$. Fix a rational point $x\in X(K)$, 
and consider the embedding $\iota : X\hookrightarrow J$ defined by $\iota(x)=0_J$.
For any field extension $\overline K\subset L$, with $L$ algebraically closed, let
$J^{\tor}\defeq J(L)^{\tor}=J(\overline K)^{\tor}$ be the torsion subgroup of $J$. The intersection $X\cap J^{\tor}$ is finite by [Raynaud]. 
Let $M$ be the cardinality of the subgroup of $J^{\tor}$ generated by $X\cap J^{\tor}$. We {\bf assume $M\ge 2$} (this $M$ will be the integer $N$ required in theorem C).

Let $p>M$ be a prime integer, $k$ a $p$-adic completion of $K$, $\overline k$ an algebraic closure of $k$, $X_k\defeq X\times_Kk$, and $X_{\overline k}\defeq X\times_K\overline k$.
Recall the exact sequence of fundamental groups (cf. $\S0$)
$$1\to \pi_1(X_{\overline k},\bar \eta)\to \pi_1(X_k,\eta)\to G_k\to 1.$$ 
Let $\Delta$ be the {\bf maximal pro-$p$ quotient} of 
$\pi_1(X_{\overline k},\bar \eta)$, and 
$$\Pi\defeq \pi_1(X_k,\eta)/\Ker (\pi_1(X_{\overline k},\bar \eta)\twoheadrightarrow \Delta )$$ 
the {\bf geometrically pro-$p$}  arithmetic fundamental group of $X_k$.

For an integer $m\ge 1$, let $\Delta_m$ be the maximal {\bf $m$-step solvable pro-$p$} quotient of 
$\pi_1(X_{\overline k},\bar \eta)$, and 
$$\Pi_m\defeq \pi_1(X_k,\eta)/\Ker (\pi_1(X_{\overline k},\bar \eta)\twoheadrightarrow \Delta_m )$$ 
the {\bf geometrically $m$-step solvable pro-$p$} arithmetic fundamental group of $X_k$. We have a commutative diagram of exact sequences
$$
\CD
@. 1 @.  1 @.\\
@. @VVV   @VVV\\
@. \Delta[m] @= \Delta [m]\\
@. @VVV  @VVV\\
1 @>>> \Delta _{m+1} @>>> \Pi _{m+1} @>>> G_k @>>> 1\\
@. @VVV  @VVV     @|\\
1 @>>> \Delta _{m} @>>> \Pi _{m} @>>> G_k @>>> 1\\
@. @VVV  @VVV\\
@. 1 @. 1\\
\endCD
\tag 2.1
$$
where $\Delta[m]\defeq \Ker (\Delta_{m+1}\twoheadrightarrow  \Delta_m)=\Ker (\Pi_{m+1}\twoheadrightarrow  \Pi_m)$ (cf. [Sa\"\i di3], $\S1$, for more details).
Further we have natural identifications 
$$\Delta=\varprojlim _{m\ge 1} \Delta_m,\ \ \ \text{and}\ \ \ \  \Pi=\varprojlim _{m\ge 1} \Pi_m.$$

Let $s_x:G_k\to \Pi$ be a section of  the projection
$\Pi\twoheadrightarrow G_k$ associated to the $k$-rational point (image in $X_k$ of) $x$, which induces sections
$s_{x,m}:G_k\to \Pi_m$ of the projections $\Pi_m\twoheadrightarrow G_k$, $\forall m\ge 1$. 
(Thus $s_x$  is defined up to conjugation by $\Delta$.)
We fix the section $s_{x,1}$ as a base point of the torsor of splittings of the exact sequence
$1 \to \Delta _{1} \to \Pi _{1} \to G_k \to 1$ ($\Pi_1$ is the geometrically abelian pro-$p$ arithmetic fundamental group of $X_k$), which is a torsor under $H^1(G_k,\Delta_1)$.





Let $y\in (X\cap J^{\tor})(K)\setminus \{0_J\}$ (the existence of $y$ follows from our assumption $M\ge 2$), $s_y:G_k\to \Pi$ a section of the projection
$\Pi\twoheadrightarrow G_k$ associated to the $k$-rational point (image in $X_k$ of) $y$, which induces sections
$s_{y,m}:G_k\to \Pi_m$ of the projections $\Pi_m\twoheadrightarrow G_k$, $\forall m\ge 1$. 
The classes $[s_{x,1}]=0$, and $[s_{y,1}]$, of the sections $s_{x,1}$, and $s_{y,1}$; respectively, in $H^1(G_k,\Delta_1)$ coincide. Indeed this follows easily from the (pro-$p$) Kummer exact sequence associated to $J$, and the fact that $\iota (y)$ is a torsion point of order prime-to-$p$ (recall $p>M$).

More generally, for $m\ge 1$, consider the following commutative diagram
$$
\CD
1@>>> \Delta [m+1] @>>> E[m+1] @>>> G_k @>>> 1\\
@. @| @VVV @Vs_{x,m}VV\\
1@>>> \Delta [m+1] @>>> \Pi_{m+1} @>>> \Pi_m @>>> 1\\
\endCD
\tag 2.2
$$
where the right square is cartesian. Thus the group extension $E[m+1]$ is the pull-back of the group extension $\Pi_{m+1}$ via the section $s_{x,m}$.

The upper exact sequence in diagram (2.2) splits. Indeed this follows from the existence of the section $s_{x,m+1}:G_k\to \Pi_{m+1}$
which lifts the section $s_{x,m}$, and induces a splitting $s_{x,m+1}:G_k\to E[m+1]$ of the group extension $E[m+1]$. We fix the section $s_{x,m+1}$ as a base point for the torsor of splittings of the group extension $E[m+1]$, which is a torsor 
under $H^1(G_k,\Delta[m+1])$; the $G_k$-module structure of $\Delta[m+1]$ is deduced from diagram (2.2). 
If $z\in X(k)$, and $s_{z,m}=s_{x,m}:G_k\to \Pi_m$, then the section $s_{z,m+1}$ gives rise to a splitting 
$s_{z,m+1}:G_k\to E[m+1]$ of the upper exact sequence in diagram (2.2), hence to a class $[s_{z,m+1}]\in H^1(G_k,\Delta[m+1])$.

Define $\Cal S_m$ to be the set of rational points $z\in X(k)$ such that $s_{x,m}(G_k)$ coincide with a decomposition group of $\Pi_m$ associated to $z$.
We have the following inclusions
$$\cdots \subseteq \Cal S_{m+1}\subseteq \Cal S_m\subseteq \cdots \subseteq\Cal S_2\subseteq \Cal S_1=X(k)\cap J^{\tor,p'}\subseteq X\cap J^{\tor}.$$
The equality $\Cal S_1=X(k)\cap J^{\tor,p'}$ follows from the (pro-$p$) Kummer exact sequence associated to $J$, and the well-known structure of $J(k)$.

\proclaim{Lemma 2.1} The equality $\bigcap _{m\ge 1}\Cal S_m=\{x\}$ holds.
\endproclaim

\demo{Proof}
Follows from [Mochizuki], Theorem C, and a limit argument using the fact that $\Pi=\varprojlim _{m\ge 1} \Pi_m$.
\qed
\enddemo

It follows from Lemma 2.1, and the above discussion, that there exists $m\ge 1$ such that 
$$\{x\}\subsetneq \Cal S_m,\ \ \  \text {and}\ \ \  \{x\}=\Cal S_{m+1}.$$
Let 
$$A\defeq \{[s_{z,m+1}] : z\in \Cal S_{m}\}\subset H^1(G_k,\Delta[m+1]).$$
Note that $\{0\}\subsetneq A$; which follows from the facts that $\{x\}\subsetneq \Cal S_m$ and $\{x\}=\Cal S_{m+1}$. Further
$\Card (A)\le \Card (\Cal S_m)\le M<p$. In particular, 
$$\exists \alpha\in H^1(G_k,\Delta[m+1])\setminus A,$$
since $H^1(G_k,\Delta[m+1])$ is $p$-primary. Thus $\alpha$ corresponds to a section $\alpha:G_k\to \Pi_{m+1}$ of the projection $\Pi_{m+1}\twoheadrightarrow G_k$, which lifts the section 
$s_{x,m}$.

\proclaim{Lemma 2.3} The section $\alpha:G_k\to \Pi_{m+1}$ is {\bf non-geometric}.
\endproclaim

\demo{Proof} Follows from the various definitions, and the fact that $\alpha\notin A$.
\qed
\enddemo

This finishes the proof of Theorem C.
\qed


$$\text{References.}$$

\noindent
[Grothendieck] Grothendieck, A., Rev\^etements \'etales et groupe fondamental, Lecture 
Notes in Math. 224, Springer, Heidelberg, 1971.

\noindent
[Hoshi] Hoshi, Y., Existence of nongeometric pro-$p$ Galois sections of hyperbolic curves, Publ. Res. Inst. Math. Sci. 46 (2010), no. 4, 829-848.

\noindent
[Mochizuki] Mochizuki, S., The Local Pro-$p$ Anabelian Geometry of Curves, Invent. Math. 138 (1999), 319-423.

\noindent
[Neukirch-Schmidt-Winberg] Neukirch, J., Schmidt, A., Winberg, K., Cohomology of Number Fields, Grundlehren der mathematischen Wissenschaften, 323, Springer, 2000.

\noindent
[Raynaud] Raynaud, M., Courbes sur une vari\'et\'e ab\'elienne et points de torsion. Invent. Math. 71, 207-233, (1983).

\noindent
[Sa\"\i di] Sa\"\i di, M., The cuspidalisation of sections of arithmetic fundamental groups, Advances in Mathematics 230
(2012) 1931-1954.

\noindent
[Sa\"\i di1]  Sa\"\i di, M., The cuspidalisation of sections of arithmetic fundamental groups II, Advances in Mathematics,
354 (2019), https://doi.org/10.1016/j.aim.2019.106737354.

\noindent
[Sa\"\i di2] Sa\"\i di, M., On the existence of non-geometric sections of arithmetic fundamental groups, 
Mathematische Zeitschrift 277, no. 1-2 (2014), 361-372.

\noindent
[Sa\"\i di3] Sa\"\i di, M., A local-global principle for torsors under geometric prosolvable fundamental groups, 
Manuscripta Mathematica 145, no. 1-2 (2014), 163-174.

\noindent
[Serre] Serre, J.-P., Cohomologie galoisienne, Seconde \'edition, 
Lecture Notes in Mathematics, 5, Springer-Verlag, Berlin-Heidelberg-New York, 1962/1963.

\noindent
[Tamagawa] Tamagawa, A., The Grothendieck conjecture for affine curves. Compositio Mathematica, 109(2) (1997), 135-194.
 
\bigskip

\noindent
Mohamed Sa\"\i di

\noindent
College of Engineering, Mathematics, and Physical Sciences

\noindent
University of Exeter

\noindent
Harrison Building

\noindent
North Park Road

\noindent
EXETER EX4 4QF

\noindent
United Kingdom

\noindent
M.Saidi\@exeter.ac.uk

\end
\enddocument